\documentclass[12pt]{amsart}%
\usepackage{amsfonts}
\usepackage{amsmath}
\usepackage{amssymb}
\usepackage{graphicx}%
\theoremstyle{plain}
\newtheorem{theorem}{Theorem}
\newtheorem*{onetheorem}{Theorem}

\newtheorem{lemma}[theorem]{Lemma}

\theoremstyle{definition}

\newtheorem{definition}{Definition}
\theoremstyle{remark}

\newtheorem*{remark}{Remark}

\begin{document}
\title[A comparison between the max and min norms]{A comparison between the max and min norms on $C^{\ast}\left(  F_{n}\right)
\otimes C^{\ast}\left(  F_{n}\right)  $}
\author{Florin R\u{a}dulescu}
\address{Department of Mathematics\\
The University of Iowa\\
Iowa City, Iowa 52242, U.S.A.}

\begin{abstract}
Let $F_{n}$, $n\geq2$, be the free group with $n$ generators, denoted by
$U_{1},U_{2},\dots,U_{n}$. Let $C^{\ast}\left(  F_{n}\right)  $ be the full
$C^{\ast}$-algebra of $F_{n}$. Let $\mathcal{X}$ be the vector subspace of the
algebraic tensor product $C^{\ast}\left(  F_{n}\right)  \otimes C^{\ast
}\left(  F_{n}\right)  $, spanned by $1\otimes1,U_{1}\otimes1,\dots
,U_{n}\otimes1,1\otimes U_{1},\dots,1\otimes U_{n}$. Let $\left\Vert
\,\cdot\,\right\Vert _{\min}$ and $\left\Vert \,\cdot\,\right\Vert _{\max}$ be
the minimal and maximal $C^{\ast}$ tensor norms on $C^{\ast}\left(
F_{n}\right)  \otimes C^{\ast}\left(  F_{n}\right)  $, and use the same
notation for the corresponding (matrix) norms induced on $M_{k}\left(
\mathbb{C}\right)  \otimes\mathcal{X}$.

Identifying $\mathcal{X}$ with the subspace of $C^{\ast}\left(  F_{2n}\right)
$ obtained by mapping $U_{1}\otimes1,\dots,1\otimes U_{n}$ into the $2n$
generators and the identity into the identity, we get a matrix norm
$\left\Vert \,\cdot\,\right\Vert _{C^{\ast}\left(  F_{2n}\right)  }$ which
dominates the $\left\Vert \,\cdot\,\right\Vert _{\max}$ norm, on $M_{k}\left(
\mathbb{C}\right)  \otimes\mathcal{X}$.

In this paper we prove that, with $N=2n+1=\dim\mathcal{X}$, we have
\[
\left\Vert X\right\Vert _{\max}\leq\left\Vert X\right\Vert _{C^{\ast}\left(
F_{2n}\right)  }\leq\left(  N^{2}-N\right)  ^{1/2}\left\Vert X\right\Vert
_{\min},\quad X\in M_{k}\left(  \mathbb{C}\right)  \otimes\mathcal{X}.
\]

\end{abstract}
\maketitle

Let $F_{n}$ be the free group on $n$ generators, $n\geq2$. Let $C^{\ast
}\left(  F_{n}\right)  $ be the full $C^{\ast}$-algebra associated with
$F_{n}$ (see, e.g., \cite{Wa}). As proved in \cite{EL}, \cite{Ta}, on the
algebraic tensor product $C^{\ast}\left(  F_{n}\right)  \otimes C^{\ast
}\left(  F_{n}\right)  $ there exist a maximal and a minimal $C^{\ast}%
$-algebra tensor norm, denoted by $\left\Vert \,\cdot\,\right\Vert _{\max}$
and $\left\Vert \,\cdot\,\right\Vert _{\min}$ respectively. Kirchberg, in
\cite{Ki}, has revived the study of the $C^{\ast}$-tensor norms on $A\otimes
A^{\operatorname{op}}$. One particular case of his very deep results shows
that the equality of the two norms on $C^{\ast}\left(  F_{\infty}\right)
\otimes C^{\ast}\left(  F_{\infty}\right)  $ is equivalent to Connes's
embedding problem (\cite{Co}).

In \cite{Pi2}, it is proven that if $E$ is a subspace of the algebraic tensor
product $A_{1}\otimes A_{2}$ of two $C^{\ast}$-algebras $A_{1}$ and $A_{2}$,
which has a basis consisting of unitaries that generate (as an algebra)
$A_{1}\otimes A_{2}$, then the complete isometry of the operator-space
structures induced on $E$ by the max and min norms implies the equality of the
$\left\Vert \,\cdot\,\right\Vert _{\max}$ and $\left\Vert \,\cdot\,\right\Vert
_{\min}$ norms on $A_{1}\otimes A_{2}$. This method is then used in \cite{Pi2}
to re-prove (and generalize) Kirchberg's theorem that $C^{\ast}\left(
F_{n}\right)  \otimes_{\max}B\left(  H\right)  =C^{\ast}\left(  F\right)
\otimes_{\min}B\left(  H\right)  $.

In this paper we consider the $N=2n+1$-dimensional subspace of $C^{\ast
}\left(  F_{n}\right)  \otimes C^{\ast}\left(  F_{n}\right)  $ generated by
$\{1\otimes1,U_{1}\otimes1,\dots,U_{n}\otimes1,1\otimes U_{1},\dots,1\otimes
U_{n}\}$. This space inherits operator-space structures (\cite{BP}, \cite{ER},
\cite{Pi1}) corresponding to the two embeddings. We denote the corresponding
norms on $\mathcal{X}\otimes M_{k}\left(  \mathbb{C}\right)  $, for all $k$ in
$\mathbb{N}$, by $\left\Vert \,\cdot\,\right\Vert _{\max}$ and $\left\Vert
\,\cdot\,\right\Vert _{\min}$.

We prove that the norm $\left\Vert \,\cdot\,\right\Vert _{\min}$ dominates the
$\left\Vert \,\cdot\,\right\Vert _{\max}$ norm, on all the tensor products in
$\mathcal{X}\otimes M_{k}\left(  \mathbb{C}\right)  $, $k\in\mathbb{N}$, by a
factor $\left(  N^{2}-N\right)  ^{1/2}$, where $N=2n+1$. More precisely, we
prove that%
\[
\left\Vert X\right\Vert _{\max}\leq\left(  N^{2}-N\right)  ^{1/2}\left\Vert
X\right\Vert _{\min},\qquad X\in M_{k}\left(  \mathbb{C}\right)
\otimes\mathcal{X}.
\]
In particular, our result, in the terminology introduced by Pisier \cite{Pi1},
also shows that the $\delta_{cb}$ (multiplicative) distance between the two
$N$-dimensional operator spaces in $C^{\ast}\left(  F_{n}\right)  \otimes
C^{\ast}\left(  F_{n}\right)  $, corresponding to the norms $\left\Vert
\,\cdot\,\right\Vert _{\max}$ and $\left\Vert \,\cdot\,\right\Vert _{\min}$,
is at most $\left(  N^{2}-N\right)  ^{1/2}$ (in general \cite{Pi1}, the
$\delta_{cb}$ distance between two finite-dimensional operator spaces of
dimension $N$ is bounded by $N$).

This work has been supported by the NSF grant DMS99-70486 and by
the Swiss National Science Foundation.

The author wishes to thank Pierre de la Harpe for the warm welcome and
mathematical discussions at the University of Gen\`{e}ve, during the summer of 2001.

\section*{Definitions}

\setcounter{equation}{-1}Let $k\in\mathbb{N}$ be a natural number and let
$\left(  e_{a,b}\right)  _{a,b=1}^{k}$ be a matrix unit in $M_{k}$ $\left(
\mathbb{C}\right)  $. Let $W_{1},W_{2},\dots,W_{2n}$ be the generators of
$F_{2n}$ and let $W_{0}=\operatorname*{Id}$. Let $\mathcal{X}$ be the subspace
of $C^{\ast}\left(  F_{2n}\right)  $ spanned by $W_{0},W_{1},\dots,W_{2n}$.
Let $X$ be an arbitrary element of $\mathcal{X}$. Then $X^{\ast}X$ has the
form%
\[
\sum_{a,b=1}^{k}\left(  \sum_{i\neq j}A_{ia,jb}W_{i}^{\ast}W_{j}^{{}}%
+B_{a,b}\operatorname*{Id}\right)  \otimes e_{a,b}.
\]

The norm $\left\Vert X\right\Vert _{C^{\ast}\left(  F_{2n}\right)  }$ for $X$
in $C^{\ast}\left(  F_{2n}\right)  $ is computed (\cite{Wa}, \cite{BP}) as the
supremum over all Hilbert spaces $H$ and all unitaries $U_{1},U_{2}%
,\dots,U_{2n}$ acting on $H$, and all $\xi=\bigoplus_{a=1}^{k}\xi_{a}$,
$\sum\left\Vert \xi_{a}\right\Vert ^{2}=1$, in $H\oplus\dots\oplus H$ ($k$
times), of the quantity%
\begin{equation}
\left\langle X^{\ast}X\xi,\xi\right\rangle =\sum_{a,b=1}^{k}\left(
\sum_{i\neq j}A_{ia,jb}\left\langle W_{i}^{\ast}W_{j}^{{}}\xi_{a},\xi
_{b}\right\rangle +B_{a,b}\left\langle \xi_{a},\xi_{b}\right\rangle \right)
.\label{eq0}%
\end{equation}
Since $C^{\ast}\left(  F_{2n}\right)  $ is residually finite \cite{Cho} (see
also \cite{Wa}, \cite{BL}), it follows that the norm of $X^{\ast}X$ might be
computed using only finite-dimensional unitaries.

Let $\tilde{V}_{1},\dots,\tilde{V}_{n}$ be the generators of a different copy
of the free group $F_{n}$. We identify $\mathcal{X}$ with a subspace of the
algebraic tensor product $C^{\ast}\left(  F_{n}\right)  \otimes C^{\ast
}\left(  F_{n}\right)  $ by mapping $1$ into $1\otimes1$, and $W_{i}$ into
$\tilde{V}_{i}\otimes1$ for $i=1,2,\dots,n$, $W_{i+n}$ into $1\otimes\tilde
{V}_{i}$ for $i=1,2,\dots,n$. With this identification, and by using again the
fact that $C^{\ast}\left(  F_{n}\right)  $ is residually finite, it follows
that the norm $\left\Vert X\right\Vert _{\max}$ viewed as an element of
$\left(  C^{\ast}\left(  F_{n}\right)  \otimes_{\max}C^{\ast}\left(
F_{n}\right)  \right)  \otimes M_{k}\left(  \mathbb{C}\right)  $ is computed
by the same supremum as the one used for $\left\Vert X\right\Vert _{C^{\ast
}\left(  F_{2n}\right)  }$, with the additional restriction, on the unitaries
$U_{1},\dots,U_{2n}$, that for $1\leq i\leq n<j\leq2n$, we have $\left[
U_{i},U_{j}\right]  =0$.

Clearly this gives (as in \cite{BP}) that $\left\Vert X\right\Vert _{C^{\ast
}\left(  F_{2n}\right)  }\geq\left\Vert X\right\Vert _{C^{\ast}\left(
F_{n}\right)  \otimes_{\max}C^{\ast}\left(  F_{n}\right)  }$. The norm
$\left\Vert X\right\Vert _{\min}$ for $X$ in $\mathcal{X}\otimes M_{k}\left(
\mathbb{C}\right)  $, viewed as an element in $C^{\ast}\left(  F_{n}\right)
\otimes_{\min}C^{\ast}\left(  F_{n}\right)  \otimes M_{k}\left(
\mathbb{C}\right)  $, is then computed by the same supremum formulas as for
$\left\Vert X\right\Vert _{\max}$, by imposing the additional condition that
the Hilbert space $H$ splits as $K_{1}\otimes K_{2}$ and there exist unitaries
$\alpha_{1},\dots,\alpha_{n}$ acting on $K_{1}$, and $\beta_{1},\dots
,\beta_{n}$ unitaries on $K_{2}$, such that $U_{i}=\alpha_{i}\otimes1$ and
$U_{i+n}=1\otimes\beta_{i}$ for $1\leq i\leq n$ (see also \cite{Vo90}).
Motivated by this we introduce the following definition:

\begin{definition}
\label{Def1}A triplet $\left(  H,\left(  U_{i}\right)  _{i=1}^{2n},\left(
\eta_{a}\right)  _{a=1}^{k}\right)  $ consisting of a Hilbert space $H$,
unitaries $\left(  U_{i}\right)  _{i=1}^{2n}$ acting on $H$ and vectors
$\left(  \eta_{a}\right)  _{a=1}^{k}$ is called \emph{in tensor position} if
there exist a Hilbert space $K$, unitaries $\tilde{U}_{1},\dots,\tilde{U}_{n}%
$, $\tilde{V}_{1},\dots,\tilde{V}_{n}$ on $K$, vectors $\left(  \tilde{\eta
}_{a}\right)  _{a=1}^{k}$ in $K\otimes K$ with the following properties.
Denote $\tilde{W}_{i}=\tilde{U}_{i}\otimes\operatorname*{Id}_{K}$ for $1\leq
i\leq n$ and $\tilde{W}_{i+n}=\operatorname*{Id}_{K}\otimes\tilde{V}_{i}$,
$1\leq i\leq n$. Also denote $U_{0}=\operatorname*{Id}_{H}$, $\tilde{W}%
_{0}=\operatorname*{Id}_{K\otimes K}$. With these notations the following
should hold true for $0\leq i,j\leq n$, $1\leq a,b\leq k$:%
\[
\left\langle U_{i}\eta_{a},U_{j}\eta_{b}\right\rangle =\left\langle \tilde
{W}_{i}\tilde{\eta}_{a},\tilde{W}_{j}\tilde{\eta}_{b}\right\rangle .
\]

\end{definition}

\section*{Main Result}

Our main result gives a comparison between the norms $\left\Vert
\,\cdot\,\right\Vert _{C^{\ast}\left(  F_{2n}\right)  }$ and $\left\Vert
\,\cdot\,\right\Vert _{\min}$ on the space $\mathcal{X}$ (and its tensor
products $\mathcal{X}\otimes M_{k}\left(  \mathbb{C}\right)  $). To do this we
use the fact that, for any triplet $\left(  H,\left(  U_{i}\right)
_{i=1}^{2n},\left(  \xi_{a}\right)  _{a=1}^{k}\right)  $, $U_{0}%
=\operatorname*{Id}$, the information contained in the matrix $\left\langle
U_{i}\xi_{a},U_{j}\xi_{b}\right\rangle $, $0\leq i,j\leq2n$, is unchanged
(except for the Gram--Schmidt matrix of $\xi_{a}$) if we replace $H$, $U_{i}$
and $\xi_{a}$ by a direct sum and linear combinations of elementary triplets
$\left(  H^{\alpha},\left(  U_{i}^{\alpha}\right)  _{i=1}^{2n},\left(  \xi
_{a}^{\alpha}\right)  _{a=1}^{k}\right)  $ having the property that the
vectors $\left\{  U_{i}^{\alpha}\xi_{a}^{\alpha}\right\}  _{i,a}$ are an
orthonormal system (with the exception of some repetitions). The following
lemma is an obvious property for triplets as in Definition \ref{Def1}:

\begin{lemma}
\label{Lem1} Let $\Lambda$ be a countable index set. Assume the triplets
\linebreak$\left(  H^{\alpha},\left(  U_{i}^{\alpha}\right)  _{i=1}%
^{2n},\left(  \eta_{a}^{\alpha}\right)  _{a=1}^{k}\right)  _{\alpha\in\Lambda
}$ are in tensor position. Let $\left(  \mu_{a}^{\alpha}\right)
_{a=1,\,\alpha\in\Lambda}^{k}$ be arbitrary complex numbers such that
$\sum_{\alpha}\left\vert \mu_{a}^{\alpha}\right\vert ^{2}\left\Vert \eta
_{a}^{\alpha}\right\Vert ^{2}<\infty$ for all $a$. Let $H=\bigoplus_{\alpha
\in\Lambda}H^{\alpha}$, let $U_{i}=\bigoplus U_{i}^{\alpha}$ and $\eta
_{a}=\bigoplus\mu_{a}^{\alpha}\eta_{a}^{\alpha}$.

Then the triplet $\left(  H,\left(  U_{i}\right)  _{i=1}^{2n},\left(  \eta
_{a}\right)  _{a=1}^{k}\right)  $ is in tensor position.
\end{lemma}

\begin{proof}
For each $\alpha\in\Lambda$, use the definition of tensor position to find a
Hilbert space $K^{\alpha}$ and unitaries $\tilde{W}_{i}^{\alpha}=\tilde{U}%
_{i}^{\alpha}\otimes\operatorname*{Id}_{K^{\alpha}}$, $1\leq i\leq n$,
$W_{i+n}^{\alpha}=\operatorname*{Id}_{K^{\alpha}}\otimes\tilde{V}_{i}^{\alpha
}$ as in Definition \ref{Def1}. Let $\tilde{K}=\bigoplus K^{\alpha}$ and
$\tilde{H}=\tilde{K}\otimes\tilde{K}\supseteq\bigoplus_{\alpha}K^{\alpha
}\otimes K^{\alpha}$. Let $\tilde{U}_{i}=\bigoplus_{\alpha}\tilde{U}%
_{i}^{\alpha}$, $\tilde{V}_{i}=\bigoplus_{\alpha}\tilde{V}_{i}^{\alpha}$ and
$\tilde{W}_{i}=\tilde{U}_{i}\otimes\operatorname*{Id}_{\tilde{K}}$, $1\leq
i\leq n$, $\tilde{W}_{i+n}=\operatorname*{Id}_{\tilde{K}}\otimes\tilde{V}_{i}%
$, $1\leq i\leq n$, $\tilde{\eta}_{a}=\bigoplus\mu_{a}^{\alpha}\tilde{\eta
}_{a}^{\alpha}$. Then the triplet $\left(  \tilde{H},\left(  \tilde{W}%
_{i}\right)  _{i=1}^{2n},\left(  \tilde{\eta}_{a}\right)  _{a=1}^{k}\right)  $
has the property that%
\[
\left\langle U_{i}\eta_{a},U_{j}\eta_{b}\right\rangle =\left\langle \tilde
{W}_{i}\tilde{\eta}_{a},\tilde{W}_{j}\tilde{\eta}_{b}\right\rangle
\]
for all $0\leq i,j\leq2n$, $a,b=1,2,\dots,k$, and hence it is in tensor position.
\end{proof}

\begin{definition}
\label{Def2}For a triplet $\left(  H,\left(  U_{i}\right)  _{i=1}^{2n},\left(
\eta_{a}\right)  _{a=1}^{k}\right)  $ \textup{(}with $U_{0}=\operatorname*{Id}%
$\textup{)}, the \emph{associated matrix} will be $X_{ia,jb}^{U}%
=X_{ia,jb}=\left\langle U_{i}\eta_{a},U_{j}\eta_{b}\right\rangle $ for $0\leq
i,j\leq2n$, $a,b=1,2,\dots,k$.
\end{definition}

Clearly $X_{ia,ib}=\left\langle \eta_{a},\eta_{b}\right\rangle $ for all $i$
and all $a,b$. Also $X_{ia,jb}=\overline{X_{jb,ia}}$ by definition.

\begin{remark}
The property in the definition of a triplet in tensor position is completely
contained in the information in the matrix $X$.

Moreover, with the notations in Lemma \textup{\ref{Lem1}}, if $X$ is the
matrix for the triplet $\left(  H,\left(  U_{i}\right)  _{i=1}^{2n},\left(
\eta_{a}\right)  _{a=1}^{k}\right)  $ and $X^{\alpha}$ is the matrix for the
triplets $\left(  H^{\alpha},\left(  U_{i}^{\alpha}\right)  _{i=1}%
^{2n},\left(  \eta_{a}^{\alpha}\right)  _{a=1}^{k}\right)  $, then we have%
\[
X_{ia,jb}=\sum_{\alpha}\mu_{a}^{\alpha}\,\overline{\mu_{b}^{\alpha}%
}\,X_{ia,jb}^{\alpha}.
\]
It is easy to construct elementary triplets in tensor position.
\end{remark}

\begin{lemma}
\label{Lem2}Let $H$ be a separable Hilbert space. Let $\varepsilon$ be a
complex number of absolute value $1$. Let $n$, $k$ be strictly positive
integers. Fix a vector $\eta$ in $H$ of length $1$. Let $\eta_{a}=\eta$ for
$a=1,\dots,k$. Let $\alpha=\left(  i_{0},j_{0}\right)  $, with $i_{0},j_{0}%
\in\left\{  0,1,\dots,2n\right\}  $, $i_{0}\neq j_{0}$. Assume $\left(
U_{i}\right)  _{i=1}^{2n}$ are unitaries such that%
\[
\bar{\varepsilon}U_{i_{0}}\eta_{a}=\bar{\varepsilon}U_{i_{0}}\eta=U_{j_{0}%
}\eta_{a}=U_{j_{0}}\eta
\]
and such that the vectors%
\[
\bar{\varepsilon}U_{i_{0}}\eta=U_{j_{0}}\eta,\ \left\{  U_{k}\eta\mid k\neq
i_{0},j_{0}\right\}
\]
are pairwise orthogonal.

Then $\left(  H,\left(  U_{i}\right)  _{i=1}^{2n},\left(  \eta_{a}\right)
_{a=1}^{k}\right)  $ is in tensor position, and the associated matrix is, for
$0\leq i,j\leq2n$, $1\leq a,b\leq k$,%
\begin{align*}
X_{ia,jb}^{\alpha,\varepsilon} &  =1\text{\qquad if }i=j,\\
X_{i_{0}a,j_{0}b}^{\alpha,\varepsilon} &  =\varepsilon,\ X_{j_{0}a,i_{0}%
b}^{\alpha,\varepsilon}=\bar{\varepsilon},\\
X_{ia,jb}^{\alpha,\varepsilon} &  =0\text{\qquad if }i\text{ or }j\text{ are
not in }\left\{  i_{0},j_{0}\right\}  \text{ and }i\neq j.
\end{align*}

\end{lemma}

\begin{proof}
It is obvious that this should be the formula for the matrix $X^{\alpha
,\varepsilon}$ associated to the triplet.

We need to construct a specific triplet in tensor position, which gives the
matrix $X^{\alpha,\varepsilon}$. To do this we split into two cases.

First we analyze the case where $0\leq i_{0}\leq n$ and $n<j_{0}\leq2n$. In
this case consider a Hilbert space $K$ of sufficiently large dimension. Let
$e_{0},e_{1},\dots$ be a basis for this Hilbert space and let $\eta$ be the
vector $e_{0}\otimes e_{0}$. With the notations from Definition \ref{Def1},
let $\tilde{W}_{i_{0}}=\operatorname*{Id}\otimes\operatorname*{Id}$,
$\tilde{W}_{j_{0}}=\bar{\varepsilon}\operatorname*{Id}\otimes
\operatorname*{Id}$ (which corresponds to the choice $\tilde{U}_{i_{0}%
}=\operatorname*{Id}$, $\tilde{V}_{j_{0}-n}=\bar{\varepsilon}%
\operatorname*{Id}$).

For $i\neq i_{0}$, $i=0,1,\dots,n$, let $\tilde{U}_{i}$ be a unitary on $K$,
such that $\left\{  \tilde{U}_{i}e_{0}\right\}  _{i\neq i_{0}}$ and $e_{0}$ is
an orthonormal system in $K$. (For example we can send $\tilde{U}_{i}e_{0}$ to
other elements in the basis.) Likewise we choose $\tilde{V}_{j}e_{0}$ such
that $\left\{  \tilde{V}_{j}e_{0}\right\}  _{j\neq j_{0}-n}$ and $e_{0}$ is an
orthonormal system. It is obvious now that the unitaries $\left(  \tilde
{U}_{i}\right)  _{i=1}^{n}$, $\left(  \tilde{V}_{j}\right)  _{j=1}^{n}$
realize a triplet in tensor position as in the statement of Lemma \ref{Lem2}.

The case $0\leq i_{0}<j_{0}\leq n$ is easier and may be treated similarly.
\end{proof}

In the next lemma we provide a decomposition of an arbitrary triplet $\left(
H,\left(  U_{i}\right)  _{i=1}^{2n},\left(  \eta_{a}\right)  _{a=1}%
^{k}\right)  $, with $H$ finite-dimensional, into elementary triplets as in
Lemma \ref{Lem2}. The drawback ot this construction is that in the
decomposition of $\left(  H,\left(  U_{i}\right)  _{i=1}^{2n},\left(  \eta
_{a}\right)  _{a=1}^{k}\right)  $, the vectors in the triplet have greater
length (by a factor of $\left(  N^{2}-N\right)  ^{1/2}$, with $N=2n+1$).

\begin{lemma}
\label{Lem3}Let $H$ be a finite-dimensional vector space. Let $n$, $k$ be
strictly positive integer numbers. Let $U_{0}=\operatorname*{Id}$,
$U_{1},\dots,U_{2n}$ be unitaries on $H$, and let $\left(  \xi_{a}\right)
_{a=1}^{k}$ be vectors in $H$.

Then there exists a triplet $\left(  \tilde{K},\left(  \tilde{U}_{i}\right)
_{i=1}^{2n},\left(  \tilde{\eta}_{a}\right)  _{a=1}^{k}\right)  $ in tensor
position, such that \textup{(}with $N=2n+1$\textup{)} we have:

\begin{enumerate}
\item \label{Lem3(1)}$\left\langle U_{i}\xi_{a},U_{j}\xi_{b}\right\rangle
=\left\langle \tilde{U}_{i}\tilde{\eta}_{a},\tilde{U}_{j}\tilde{\eta}%
_{b}\right\rangle $, $i\neq j$,

\item \label{Lem3(2)}$\left\langle \tilde{\eta}_{a},\tilde{\eta}%
_{b}\right\rangle =\left(  N^{2}-N\right)  \left\langle \xi_{a},\xi
_{b}\right\rangle $
\end{enumerate}

\noindent for all $a,b=1,2,\dots,k$ and for all $i,j=0,1,\dots,2n$
\textup{(}and $i\neq j$\textup{)}.
\end{lemma}

\begin{proof}
Let $\left(  e_{t}\right)  _{t\in T}$ be an orthonormal basis for $H$ and let
$\lambda_{i,a}^{t}$ be the components of the vector $U_{i}\xi_{a}$ in this
basis for $i=0,1,\dots,2n$, $a=1,\dots,k$, $t\in T$. Then we have that%
\begin{equation}
\left\langle U_{i}\xi_{a},U_{j}\xi_{b}\right\rangle =\sum_{t}\lambda_{i,a}%
^{t}\,\overline{\lambda_{j,b}^{t}}\,,\qquad i,j=0,1,\dots,2n,\ a,b=1,\dots
,k.\label{eq1}%
\end{equation}
The usual factorization formula \cite{Pe} gives, with $\varepsilon=\sqrt{-1}$,
that for all $i,j=0,1,\dots,2n$ and for all $a,b=1,\dots,k$ we have that%
\begin{equation}
\lambda_{i,a}^{t}\,\overline{\lambda_{j,b}^{t}}=\frac{1}{4}\sum_{s=0}%
^{3}\varepsilon^{s}\left(  \lambda_{i,a}^{t}+\varepsilon^{s}\lambda_{j,a}%
^{t}\right)  \,\overline{\left(  \lambda_{i,b}^{t}+\varepsilon^{s}%
\lambda_{j,b}^{t}\right)  }\,.\label{eq2}%
\end{equation}
Note also that the following holds:%
\begin{equation}
\frac{1}{4}\sum_{s=0}^{3}\left(  \lambda_{i,a}^{t}+\varepsilon^{s}%
\lambda_{j,a}^{t}\right)  \,\overline{\left(  \lambda_{i,b}^{t}+\varepsilon
^{s}\lambda_{j,b}^{t}\right)  }=\lambda_{i,a}^{t}\,\overline{\lambda_{i,b}%
^{t}}+\lambda_{j,a}^{t}\,\overline{\lambda_{j,b}^{t}}\,.\label{eq3}%
\end{equation}
For a given pair $\alpha=\left(  i,j\right)  $, $0\leq i<j\leq n$,
$a,b=1,\dots,k$, $t\in T$, and $s=0,1,2,3$, we let%
\[
\theta_{\alpha,a}^{t,s}=\lambda_{i,a}^{t}+\varepsilon^{s}\lambda_{j,a}^{t}.
\]
With these notations the relations (\ref{eq2}) and (\ref{eq3}) become
respectively%
\begin{align}
\left\langle U_{i}\xi_{a},U_{j}\xi_{b}\right\rangle  &  =\sum_{t}\lambda
_{i,a}^{t}\lambda_{j,b}^{t}\label{eq4}\\
&  =\sum_{t,s}\varepsilon^{s}\theta_{\alpha,a}^{t,s}\,\overline{\theta
_{\alpha,b}^{t,s}}\,,\nonumber\\
\sum_{t,s}\theta_{\alpha,a}^{t,s}\,\overline{\theta_{\alpha,b}^{t,s}} &
=\sum_{t}\lambda_{i,a}^{t}\,\overline{\lambda_{i,b}^{t}}+\sum_{t}\lambda
_{j,a}^{t}\,\overline{\lambda_{j,b}^{t}}\label{eq5}\\
&  =\left\langle U_{i}\xi_{a},U_{i}\xi_{b}\right\rangle +\left\langle U_{j}%
\xi_{a},U_{j}\xi_{b}\right\rangle =2\left\langle \xi_{a},\xi_{b}\right\rangle
.\nonumber
\end{align}
The relations (\ref{eq4}) and (\ref{eq5}) hold for all $0\leq i<j\leq2n$, and
all $a,b=1,2,\dots,k$.

For each fixed $t\in T$, $\alpha=\left(  i_{0},j_{0}\right)  $, $0\leq
i_{0}<j_{0}\leq2n$, and each $s=0,1,2,3$, let $\left(  H^{\alpha,s,t},\left(
U_{i}^{\alpha,s,t}\right)  _{i=1}^{2n},\left(  \eta_{a}^{\alpha,s,t}\right)
_{a=1}^{k}\right)  $ be the triplet constructed in Lemma \ref{Lem2} for
$\varepsilon=\varepsilon^{s}$. (This triplet does not depend on $t$, but for
each $t$ we consider one copy.) The matrix associated to this triplet is
defined by
\begin{align}
X_{ia,jb}^{\alpha,s,t} &  =0\text{\qquad if }\left\{  i,j\right\}
\nsubseteq\left\{  i_{0},j_{0}\right\}  \text{ and }i\neq j,\label{eq6}\\
X_{ia,ib}^{\alpha,s,t} &  =1,\nonumber\\
X_{i_{0}a,j_{0}b}^{\alpha,s,t} &  =\varepsilon^{s},\ X_{j_{0}a,i_{0}b}%
^{\alpha,s,t}=\overline{\varepsilon^{s}}\nonumber\\
&  \text{\qquad for all }a,b=1,2,\dots,k.\nonumber
\end{align}
Let $\Lambda$ be the set of pairs%
\[
\Lambda=\left\{  \left(  i,j\right)  \mid0\leq i<j\leq2n\right\}  .
\]
Let $\mu_{a}^{\alpha,s,t}=\theta_{\alpha,a}^{s,t}$ for all $\alpha\in\Lambda$,
$s=0,1,2,3$, $t\in T$. We apply Lemma \ref{Lem1} (and the following remark) to
the direct sum of the triplets $\left(  H^{\alpha,s,t},\left(  U_{i}%
^{\alpha,s,t}\right)  _{i=1}^{2n},\left(  \eta_{a}^{\alpha,s,t}\right)
_{a=1}^{k}\right)  $. In the direct sum $\tilde{H}=\bigoplus_{\alpha
,s,t}H^{\alpha,s,t}$, $\tilde{U}_{i}=\bigoplus_{\alpha,s,t}U_{i}^{\alpha,s,t}%
$, $i=1,2,\dots,2n$, we consider the vectors $\tilde{\eta}_{a}=\bigoplus
_{\alpha,s,t}\mu_{a}^{\alpha,s,t}\eta_{a}^{\alpha,s,t}$.

By Lemma \ref{Lem1}, for fixed $i_{0}<j_{0}$, $a,b=1,2,\dots,k$, we have%
\[
\left\langle \tilde{U}_{i_{0}}\tilde{\eta}_{a},U_{j_{0}}\tilde{\eta}%
_{b}\right\rangle =\sum_{\alpha,s,t}\mu_{a}^{\alpha,s,t}\mu_{b}^{\alpha
,s,t}X_{i_{0}a,j_{0}b}^{\alpha,s,t}.
\]
By the relation (\ref{eq6}), and since $i_{0}<j_{0}$, an entry in the matrix
$X_{i_{0}a,j_{0}b}^{\alpha,s,t}$ is nonzero only when $\alpha$ is equal to
$\left(  i_{0},j_{0}\right)  $, and is equal in this case to $\varepsilon^{s}%
$. Thus, with $\alpha_{0}=\left(  i_{0},j_{0}\right)  $ and using the relation
(\ref{eq4}), we obtain%
\begin{align}
\left\langle \tilde{U}_{i_{0}}\eta_{a},\tilde{U}_{j_{0}}\eta_{b}\right\rangle
&  =\sum_{s,t}\varepsilon^{s}\mu_{a}^{\alpha_{0},s,t}\mu_{b}^{\alpha_{0}%
,s,t}\label{eq7}\\
&  =\sum_{s,t}\varepsilon^{s}\theta_{\alpha_{0},a}^{t,s}\,\overline
{\theta_{\alpha_{0},b}^{t,s}}\nonumber\\
&  =\left\langle U_{i_{0}}\xi_{a},U_{j_{0}}\xi_{b}\right\rangle \text{\qquad
for all }a,b=1,\dots,k.\nonumber
\end{align}

Since also $\left\langle U_{j_{0}}\xi_{b},U_{i_{0}}\xi_{a}\right\rangle
=\overline{\left\langle U_{i_{0}}\xi_{a},U_{j_{0}}\xi_{b}\right\rangle }$ and
similarly for $\tilde{U}_{i}\tilde{\eta}_{a}$, it follows that relation
(\ref{eq7}) holds for all $i_{0}\neq j_{0}$, $0\leq i_{0},j_{0}\leq2n$.

Similar computations yield the value of $\left\langle \tilde{\eta}_{a}%
,\tilde{\eta}_{b}\right\rangle $. Indeed, by the relation (\ref{eq5}) we have%
\begin{align*}
\left\langle \tilde{\eta}_{a},\tilde{\eta}_{b}\right\rangle  &  =\sum
_{\alpha,s,t}\mu_{a}^{\alpha,s,t}\,\overline{\mu_{b}^{\alpha,s,t}}\\
&  =\sum_{\alpha\in\Lambda}\sum_{s,t}\theta_{\alpha,a}^{t,s}\,\overline
{\theta_{\alpha,b}^{t,s}}\\
&  =\sum_{\alpha\in\Lambda}2\left\langle \xi_{a},\xi_{b}\right\rangle
=\frac{N^{2}-N}{2}\cdot2\left\langle \xi_{a},\xi_{b}\right\rangle =\left(
N^{2}-N\right)  \left\langle \xi_{a},\xi_{b}\right\rangle .
\end{align*}
By Lemmas \ref{Lem1} and \ref{Lem2}, the triplet $\left(  \tilde{H},\left(
\tilde{U}_{i}\right)  _{i=1}^{2n},\left(  \tilde{\eta}_{a}\right)  _{a=1}%
^{k}\right)  $ is in tensor position. This completes the proof of Lemma
\ref{Lem3}.
\end{proof}

We now can prove the main result. We will show that on $\mathcal{X}%
=\operatorname*{Sp}\{1\otimes1,U_{1}\otimes1,\dots,U_{n}\otimes1,1\otimes
U_{1},\dots,1\otimes U_{n}\}$, the matrix norm structures induced by the norms
$\left\Vert \,\cdot\,\right\Vert _{\max}$ and $\left\Vert \,\cdot\,\right\Vert
_{\min}$ on $C^{\ast}\left(  F_{n}\right)  \otimes C^{\ast}\left(
F_{n}\right)  $ are comparable by a factor $\left(  N^{2}-N\right)  ^{1/2}$.

In particular this shows (in the terminology introduced in \cite{Pi1}) that
the $\delta_{cb}$ multiplicative distance between the two operator spaces is
less than $\left(  N^{2}-N\right)  ^{1/2}$. (By \cite{Pi1}, this distance is
at most $N$.)

\begin{onetheorem}
Let $n$, $k$ be integers, $n\geq2$, $k\geq1$. Let $F_{n}$ be the free group on
$n$ generators $V_{1},V_{2},\dots,V_{n}$. Consider the vector subspace
$\mathcal{X}$ of $C^{\ast}\left(  F_{n}\right)  \otimes C^{\ast}\left(
F_{n}\right)  $ spanned by $\{1\otimes1,V_{1}\otimes1,\dots,V_{n}%
\otimes1,1\otimes V_{1},\dots,1\otimes V_{n}\}$. Clearly $\mathcal{X}$ has
dimension $N=2n+1$.

By embedding $\mathcal{X}$ into $C^{\ast}\left(  F_{n}\right)  \otimes_{\min
}C^{\ast}\left(  F_{n}\right)  $ or $C^{\ast}\left(  F_{n}\right)
\otimes_{\max}C^{\ast}\left(  F_{n}\right)  $ respectively, we get two
corresponding norms on $\mathcal{X}\otimes M_{k}\left(  \mathbb{C}\right)  $,
denoted by $\left\Vert \,\cdot\,\right\Vert _{\max}$ and $\left\Vert
\,\cdot\,\right\Vert _{\min}$.

Let $F_{2n}$ be the free group with $2n$ generators $W_{1},\dots,W_{2n}$. We
also identify $\mathcal{X}$ with a subspace of the full $C^{\ast}$-algebra
$C^{\ast}\left(  F_{2n}\right)  $ by mapping $1\otimes1$ into $1$ and
$V_{1}\otimes1,\dots,V_{n}\otimes1$ into $W_{1},\dots,W_{n}$ and $1\otimes
V_{1},\dots,1\otimes V_{n}$ into $W_{n+1},\dots,W_{2n}$ respectively. For $X$
in $\mathcal{X}\otimes M_{k}\left(  \mathbb{C}\right)  $ we denote the
corresponding norm coming from this embedding by $\left\Vert X\right\Vert
_{C^{\ast}\left(  F_{2n}\right)  }$.

Then for all $X$ in $\mathcal{X}\otimes M_{k}\left(  \mathbb{C}\right)  $ we
have%
\[
\left\Vert X\right\Vert _{\min}\leq\left\Vert X\right\Vert _{\max}%
\leq\left\Vert X\right\Vert _{C^{\ast}\left(  F_{2n}\right)  }\leq\left(
N^{2}-N\right)  ^{1/2}\left\Vert X\right\Vert _{\min}.
\]

\end{onetheorem}

\begin{proof}
Let $\left(  e_{a,b}\right)  _{a,b=1}^{k}$ be a matrix unit in $M_{k}\left(
\mathbb{C}\right)  $ and let%
\begin{equation}
X=\sum_{r,s=1}^{k}\sum_{i=0}^{2n}\lambda_{r,s}^{i}W_{i}\otimes e_{r,s}%
,\qquad\lambda_{r,s}^{i}\in\mathbb{C}\text{,}\label{eq8}%
\end{equation}
be an arbitrary element in $M_{k}\left(  \mathbb{C}\right)  \otimes\mathbb{C}%
$. (We denote by $W_{0}$ the identity.)

Then obviously%
\begin{equation}
X^{\ast}X=\sum_{a,b=1}^{k}\left(  \sum_{\substack{i,j=0\\i\neq j}%
}^{2n}A_{ia,jb}W_{i}^{\ast}W_{j}^{{}}+B_{a,b}\operatorname*{Id}\right)
\otimes e_{a,b},\label{eq9}%
\end{equation}
where for $i\neq j$, $i,j=0,\dots,2n$, $1\leq a,b\leq k$, we have%
\begin{align}
A_{ia,jb} &  =\sum_{r=1}^{k}\lambda_{r,a}^{i}\,\overline{\lambda_{r,b}^{i}%
}\,,\label{eq10}\\
B_{a,b} &  =\sum_{r=1}^{k}\sum_{i=0}^{2n}\lambda_{r,a}^{i}\,\overline
{\lambda_{r,b}^{i}}\,.\label{eq11}%
\end{align}
Clearly the matrix $\sum_{a,b}B_{a,b}\otimes e_{a,b}$ is positive. By
definition, the $C^{\ast}\left(  F_{2n}\right)  $-norm of a noncommutative
polynomial $P$ in $\operatorname*{Id},W_{1},\dots,W_{2n}$ is computed by
taking the supremum, over all unitaries $U_{1},\dots,U_{2n}$, of the norms of
the operators obtained by replacing in $P$ the unitaries $W_{i}$ by $U_{i}$,
$i=1,2,\dots,2n$.

By \cite{Cho}, $C^{\ast}\left(  F_{2n}\right)  $ is residually finite
(\cite{Wa}, \cite{Vo90}, \cite{BL}), and hence we can restrict to a supremum
over unitaries acting on finite-dimensional vector spaces.

As a consequence, the square of the $\left\Vert \,\cdot\,\right\Vert
_{C^{\ast}\left(  F_{2n}\right)  }$ norm of the element $X$ is computed as the
supremum, over all finite-dimensional Hilbert spaces $H$, all $2n$-tuples of
unitaries $U_{1},\dots,U_{2n}$ acting on $H$ and all vectors $\xi=\left(
\xi_{a}\right)  _{a=1}^{k}$ in $H\oplus H\oplus\dots\oplus H$, $\sum_{a=1}%
^{k}\left\Vert \xi_{a}\right\Vert ^{2}=1$, of the quantities%
\begin{align*}
\left\Vert X\xi\right\Vert ^{2} &  =\left\langle X^{\ast}X\xi,\xi\right\rangle
\\
&  =\sum_{a,b}\left(  \sum_{\substack{i,j=0\\i\neq j}}^{2n}A_{ia,jb}%
\left\langle U_{j}\xi_{a},U_{i}\xi_{b}\right\rangle +B_{a,b}\left\langle
\xi_{a},\xi_{b}\right\rangle \right)  .
\end{align*}
Similarly the norm $\left\Vert X\right\Vert _{\min}$ will be computed as the
supremum of the same quantities, with the additional condition that the
unitaries $U_{1},\dots,U_{2n}$ are represented on a Hilbert space
$H=K_{1}\otimes K_{2}$, and there are unitaries $\alpha_{1},\dots,\alpha_{n}$,
respectively $\beta_{1},\dots,\beta_{n}$, on $K_{1}$, respectively $K_{2}$,
such that $U_{i}=\alpha_{i}\otimes1$, $U_{i+n}=1\otimes\beta_{i}$, $1\leq
i\leq n$.

Hence for every $\varepsilon>0$, there exists a triplet $\left(  H,\left(
U_{i}\right)  _{i=1}^{2n},\left(  \xi_{a}\right)  _{a=1}^{k}\right)  $
consisting of a finite-dimensional vector space, $2n$ unitaries on $H$ and $k$
vectors in $H$, such that (with $U_{0}=\operatorname*{Id}$)%
\begin{multline}
\left\Vert X^{\ast}X\right\Vert _{C^{\ast}\left(  F_{2n}\right)  }%
-\varepsilon\label{eq12}\\
\leq\sum_{a,b=1}^{k}\left(  \sum_{\substack{i,j=0\\i\neq j}}^{2n}%
A_{ia,jb}\left\langle U_{j}\xi_{a},U_{i}\xi_{b}\right\rangle +B_{a,b}%
\left\langle \xi_{a},\xi_{b}\right\rangle \right)  .
\end{multline}
By Lemma \ref{Lem3} we can find a triplet in tensor position, $\left(
\tilde{H},\left(  \tilde{U}_{i}\right)  _{i=1}^{2n},\left(  \tilde{\eta}%
_{a}\right)  _{a=1}^{k}\right)  $, consisting of unitaries $\tilde{U}_{i}$ on
$\tilde{H}$ (with $\tilde{U}_{0}=\operatorname*{Id}$) and vectors $\tilde
{\eta}_{a}\in\tilde{H}$ such that for all $a$, $b$,%
\begin{align}
\left\langle U_{j}\xi_{a},U_{i}\xi_{b}\right\rangle  &  =\left\langle
\tilde{U}_{j}\tilde{\eta}_{a},\tilde{U}_{i}\tilde{\eta}_{b}\right\rangle
,\qquad\qquad i\neq j,\ i,j=0,\dots,2n,\label{eq13}\\
\left\langle \tilde{\eta}_{a},\tilde{\eta}_{b}\right\rangle  &  =\left(
N^{2}-N\right)  ^{\frac{1}{2}}\left\langle \xi_{a},\xi_{b}\right\rangle
.\label{eq14}%
\end{align}
The relation (\ref{eq14}) implies that%
\[
\sum_{a}\left\Vert \tilde{\eta}_{a}\right\Vert ^{2}=\sum_{a}\left\langle
\tilde{\eta}_{a},\tilde{\eta}_{a}\right\rangle =\left(  N^{2}-N\right)
\sum_{a}\left\Vert \xi_{a}\right\Vert ^{2}=\left(  N^{2}-N\right)  .
\]
Thus, by the definition of the norm $\left\Vert X\right\Vert _{\min}$, and
since $\left(  \tilde{U}_{i}\right)  _{i=1}^{2n}$ are in tensor position, it
follows that%
\begin{multline}
\sum_{a,b}^{k}\left(  \sum_{\substack{i,j=0\\i\neq j}}^{2n}A_{ia,jb}%
\left\langle \tilde{U}_{j}\tilde{\eta}_{a},\tilde{U}_{i}\tilde{\eta}%
_{b}\right\rangle +B_{a,b}\left\langle \tilde{\eta}_{a},\tilde{\eta}%
_{b}\right\rangle \right)  \label{eq15}\\
\leq\left(  N^{2}-N\right)  \left\Vert X\right\Vert _{\min}^{2}.
\end{multline}
Moreover, the relation (\ref{eq14}) and the fact that the matrix $\sum
_{a,b}B_{a,b}\otimes e_{a,b}$ is positive imply that the right-hand side in
the inequality (\ref{eq12}) is less than the left-hand side in the inequality
in (\ref{eq15}). Hence%
\[
\left\Vert X^{\ast}X\right\Vert _{C^{\ast}\left(  F_{2n}\right)  }%
-\varepsilon\leq\left(  N^{2}-N\right)  \left\Vert X\right\Vert _{\min}^{2}.
\]
Since $\varepsilon$ is arbitrary, the result follows.
\end{proof}

\end{document}